\documentclass[12pt]{article}
\usepackage{amssymb}
\usepackage{latexsym}
\usepackage{amsmath,amssymb,latexsym}

\newfont{\blb}{msbm10 scaled\magstep1}
\newtheorem{theo}{Theorem}[section]

\newtheorem{prop}[theo]{Proposition}
\newtheorem{lemm}[theo]{Lemma}
\newtheorem{coro}[theo]{Corollary}
\title{A left $3$-Engel element \\ whose normal closure is not nilpotent}
\date{}
\author{Marialaura Noce\footnote{University of Salerno - University of the Basque Country, \texttt{mnoce@unisa.it}}, Gareth Tracey\footnote{University of Bath, \texttt{g.m.tracey@bath.ac.uk }}, and Gunnar Traustason\footnote{University of Bath, \texttt{g.traustason@bath.ac.uk}}}

\begin{document}
	
\maketitle
\begin{abstract}
\noindent We give an example of a locally nilpotent group $G$ containing a left $3$-Engel element $x$ where
$\langle x\rangle^{G}$ is not nilpotent. 
\end{abstract}

\vspace{0.16cm}
\noindent{\footnotesize 2010 Mathematics Subject Classification:  20F45, 20F40 }
\vspace{0.16cm}

\section{Introduction}

\mbox{}\\
Let $G$ be a group. An element $a\in G$ is a left Engel element in $G$, if for
each $x\in G$ there exists a non-negative integer $n(x)$ such that
      $$[[[x,\underbrace {a],a],\ldots ,a]}_{n(x)}=1.$$
If $n(x)$ is bounded above by $n$ then we say that $a$ is a left $n$-Engel element
in $G$. It is straightforward to see that any element of the Hirsch-Plotkin
radical $HP(G)$ of $G$ is a left Engel element and the converse is known
to be true for some classes of groups, including solvable groups and 
finite groups (more generally groups satisfying the maximal condition on
subgroups) [3,1]. The converse is however not true in general and this is the case
even for bounded left Engel elements. In fact whereas one sees readily that 
a left $2$-Engel element is always in the Hirsch-Plotkin radical this
is still an open question for left  $3$-Engel elements. Recently there has been a breakthrough and in [6] it is shown that any left $3$-Engel element of odd order is contained in $HP(G)$. From [10] one also knows that in order to generalise this to left $3$-Engel elements of any finite order it suffices to deal with elements of
order $2$. \\ \\
It was observed by William
Burnside [2] that every element in a group of exponent $3$  is a left $2$-Engel
element and so the fact that every left $2$-Engel element lies in the Hirsch-Plotkin radical can be seen as the underlying reason why groups of exponent
$3$ are locally finite. For groups of $2$-power exponent there is a close link
with left Engel elements. If  $G$ is a group of exponent 
$2^{n}$ then it is not difficult to see that any element $a$ in $G$ of order $2$ is a left $(n+1)$-Engel element of $G$ (see the introduction of [11] for details). For sufficiently large $n$ we know that the variety of groups of exponent $2^{n}$ is not locally finite [5,7]. As a result one can see (for example in [11]) that it follows that for sufficiently large $n$ we do not have in general that a left $n$-Engel element is contained in the Hirsch-Plotkin radical. Using the fact that groups of exponent $4$ are locally finite [9], one can also see that if all left $4$-Engel elements of a group $G$ of exponent $8$ are in
$HP(G)$ then $G$ is locally finite. \\ \\
Swapping the role of $a$ and $x$ in the definition of a left Engel element we get the notion of a right Engel element. Thus an element $a\in G$ is a right Engel element, if for each $x\in G$ there exists a non-negative 
integer $n(x)$ such that 
    $$[[[a,\underbrace {x],x],\ldots ,x]}_{n(x)}=1.$$
If $n(x)$ is bounded above by $n$, we say that $a$ is a right $n$-Engel element. By a classical result
of Heineken [4] one knows that if $a$ is a right $n$-Engel element in $G$ then $a^{-1}$ is a left $(n+1)$-Engel
element. \\ \\
In [8] M. Newell proved that if $a$ is a right $3$-Engel element in $G$ then $a\in HP(G)$ and in
fact he proved the stronger result that $\langle a\rangle^{G}$ is nilpotent of class at most $3$. The natural
question arises whether the analogous result holds for left $3$-Engel elements. In this paper we show this is
not the case by giving an example of a locally finite $2$-group with a left $3$-Engel element $a$ such that
$\langle a\rangle^{G}$ is not nilpotent. 
\section{An example of a left $3$-Engel element $x$  where $\langle x\rangle^{G}$ is not nilpotent}
Our construction will be based on an example of a Lie algebra given in [12].
Let  ${\mathbb F}$  be the field of order $2$ and consider a $4$-dimensional vector space $V={\mathbb F}x+{\mathbb F}u+{\mathbb F}v+
{\mathbb F}w$ where
             $$u\cdot v=u,\ v\cdot w=w,\ w\cdot u=v,\ u\cdot x=0,\, v\cdot x=0,\, w\cdot x=u.$$
We then extend the product linearly on $V$. One can check that $V$ is a Lie algebra with a trivial center and where $W={\mathbb F}u+{\mathbb F}v+{\mathbb F}w$ is a simple ideal. \\ \\
Let $E=\langle \mbox{ad}(x),\mbox{ad}(u),\mbox{ad}(v),\mbox{ad}(w)\rangle\leq \mbox{End}(V)$ be the associative enveloping algebra of $V$. Let 
$$\begin{array}{llll}
     e_{1}=\mbox{ad}(w), & e_{2}=\mbox{ad}(w)^{2}, & e_{3}=\mbox{ad}(w)^{3},  & e_{4}=\mbox{ad}(v), \\
   e_{5}=\mbox{ad}(v)\mbox{ad}(w), & e_{6}=\mbox{ad}(v)\mbox{ad}(w)^{2},  &  e_{7}=\mbox{ad}(u), & e_{8}=\mbox{ad}(u)\mbox{ad}(w), \\
   e_{9}=\mbox{ad}(u)\mbox{ad}(w)^{2}, & e_{10}=\mbox{ad}(x)\mbox{ad}(v), & e_{11}=\mbox{ad}(x)\mbox{ad}(w), & e_{12}=\mbox{ad}(x)\mbox{ad}(w)^{2}.
\end{array}$$
\begin{lemm} The associative  enveloping algebra $E$ is $12$-dimensional with basis $e_{1},\ldots ,e_{12}$.
\end{lemm}
{\bf Proof}\ \ We first show that $E$ is spanned by products of the form $$\mbox{ad}(x)^{\epsilon}\cdot\mbox{ad}(u)^{r}\cdot\mbox{ad}(v)^{s}\cdot\mbox{ad}(w)^{t}$$
where  $\epsilon, r,s,t$ are non-negative integers. To see this we need to show that any product $\mbox{ad}(y_{1})\cdots \mbox{ad}(y_{m})$, with 
$y_{1},\ldots ,y_{m}\in \{x,u,v,w\}$ can be written as a linear combination of elements of the required form. We use induction on $m$. This is obvious when
$m=1$. Now suppose $m\geq 2$ and that the statement is true for all shorter products. Suppose there are $\epsilon$ entries of $x$, $r$ entries of
$u$, $s$ entries of $v$ and $t$ entries of $w$ in the product $\mbox{ad}(y_{1})\cdots \mbox{ad}(y_{m})$. Using the fact that $\mbox{ad}(y_{i})\mbox{ad}(y_{j})=\mbox{ad}(y_{j})\mbox{ad}(y_{i})+\mbox{ad}(y_{i}y_{j})$, we see that modulo shorter products we have 
                $$\mbox{ad}(y_{1})\cdots \mbox{ad}(y_{m})=\mbox{ad}(x)^{\epsilon}\mbox{ad}(u)^{r}\mbox{ad}(v)^{s}\mbox{ad}(w)^{t}.$$
Hence the statement is true for products of length $m$. This finishes the inductive proof of our claim.  \\ \\
From the fact that $\mbox{ad}(x)^{2}=\mbox{ad}(x)\mbox{ad}(u)=0$, $\mbox{ad}(v)^{2}=\mbox{ad}(v)$, $\mbox{ad}(u)^{2}=\mbox{ad}(x)\mbox{ad}(v)$ and $\mbox{ad}(u)\mbox{ad}(v)=\mbox{ad}(u)+\mbox{ad}(x)\mbox{ad}(w)$, we can assume that $0\leq \epsilon,r,s\leq 1$ and
that if $\mbox{ad}(u)$ is included then we can assume that neither $\mbox{ad}(x)$ nor $\mbox{ad}(v)$ is included. This together with $\mbox{ad}(x)=\mbox{ad}(x)\mbox{ad}(v)$ and $\mbox{ad}(w)^{4}=\mbox{ad}(v)\mbox{ad}(w)^{3}=\mbox{ad}(u)\mbox{ad}(w)^{3}=
\mbox{ad}(x)\mbox{ad}(w)^{3}=0$ shows that $E$ is generated by $e_{1},\ldots ,e_{12}$. It remains to see that these elements are 
linearly independent. Suppose $\alpha_{1}e_{1}+\cdots +\alpha_{12}e_{12}=0$ for some $\alpha_{1},\ldots ,\alpha_{12}\in {\mathbb F}$. Then 
     $$0=x(\alpha_{1}e_{1}+\ldots \alpha_{12}e_{12})=\alpha_{1}u+\alpha_{2}v+\alpha_{3}w$$
gives that $\alpha_{1}=\alpha_{2}=\alpha_{3}=0$. Then 
   $$0=u(\alpha_{4}e_{4}+\ldots +\alpha_{12}e_{12})=\alpha_{4}u+\alpha_{5}v+\alpha_{6}w$$
implies that $\alpha_{4}=\alpha_{5}=\alpha_{6}=0$. Likewise 
     $$0=v(\alpha_{7}e_{7}+\ldots +\alpha_{12}e_{12})=\alpha_{7}u+\alpha_{8}v+\alpha_{9}w$$
giving $\alpha_{7}=\alpha_{8}=\alpha_{9}=0$. Finally 
    $$0=w(\alpha_{10}e_{10}+\alpha_{12}e_{11}+\alpha_{12}e_{12})=\alpha_{10}u+\alpha_{11}v+\alpha_{12}w$$
and thus $\alpha_{10}=\alpha_{11}=\alpha_{12}=0$. $\Box$ 
\\ \\
We use this example to construct a certain locally nilpotent Lie algebra over ${\mathbb F}$ of countably infinite dimension. For ease of notation it will be useful 
to introduce the following modified union of subsets of ${\mathbb N}$.  We let 
           $$A\sqcup B=\left[ \begin{array}{l}
                                     A\cup B\mbox{\ (if }A\cap B=\emptyset) \\
                                                 \emptyset\mbox{\ (otherwise)}
                                        \end{array}\right.$$
For each non-empty subset $A$ of ${\mathbb N}$ we let $W_{A}$ be a copy of the vector space $W$. That is $W_{A}=\{z_{A}:z\in W\}$ with
addition $z_{A}+t_{A}=(z+t)_{A}$. We then take the direct sum of these
     $$W^{*}=\bigoplus_{\emptyset\not =A\subseteq {\mathbb N}}W_{A}$$
that we turn into a Lie algebra with multiplication 
         $$z_{A}\cdot t_{B}=(zt)_{A\sqcup B}$$
when $z_{A}\in W_{A}$ and $t_{B}\in W_{B}$ that is then extended linearly on $W^{*}$. The interpretation here is that $z_{\emptyset}=0$. Finally we extend this to a semidirect product with ${\mathbb F}x$ 
             $$V^{*}=W^{*}\oplus {\mathbb F}x$$
induced from the action $z_{A}\cdot x=(zx)_{A}$.  \\ \\
Notice that $V^{*}$ has basis 
       $$\{x\}\cup\{u_{A},v_{A},w_{A}:\emptyset\not =A\subseteq {\mathbb N}\}$$
and that 
        $$u_{A}\cdot u_{B}=v_{A}\cdot v_{B}=w_{A}\cdot w_{B}=0,$$
     $$u_{A}\cdot x=0,\ v_{A}\cdot x=0,\ w_{A}\cdot x=u_{A}$$
and
    $$u_{A}\cdot v_{B}=u_{A\sqcup B},\ v_{A}\cdot w_{B}=w_{A\sqcup B},\ w_{A}\cdot u_{B}=v_{A\sqcup B}.$$
Notice that any finitely generated subalgebra of $V^{*}$ is contained in some $S=\langle x, u_{A_{1}}, \ldots ,
u_{A_{r}}, v_{B_{1}},\ldots ,v_{B_{s}},w_{C_{1}},\ldots ,w_{C_{t}}\rangle$. From the fact that $zxx=0$ for all $z\in V^{*}$ it follows that $S$ is nilpotent
of class at most $2(r+s+t)$. Hence $V^{*}$ is locally nilpotent. The next aim is to find a group $G\leq \mbox{GL}(V^{*})$ containing $1+\mbox{ad}(x)$ where
$1+\mbox{ad}(x)$ is a left $3$-Engel element in $G$ but where $\langle 1+\mbox{ad}(x)\rangle^{G}$ is not nilpotent. The next lemma is a preparation for this. 
\begin{lemm} The adjoint linear operator $\mbox{ad}(x)$ on $V^{*}$ satisfies:\\
(a) $\mbox{ad}(x)^{2}=0$. \\
(b) $\mbox{ad}(x)\mbox{ad}(y)\mbox{ad}(x)=0$ for all $y\in V^{*}$.
\end{lemm}
{\bf Proof}\ \ (a) Follows from the fact that $x\cdot x=u_{A}\cdot x=v_{A}\cdot x=0$ and $(w_{A}\cdot x)\cdot
x=u_{A}\cdot x=0$. \\ \\
(b) Follows from $w_{A}\cdot x\cdot u_{B}=u_{A}\cdot u_{B}=0$, $w_{A}\cdot x\cdot v_{B}\cdot x=u_{A}\cdot v_{B}\cdot x =u_{A\sqcup B}\cdot x=0$ and $w_{A}\cdot x\cdot w_{B}\cdot x=u_{A}\cdot w_{B}\cdot x=v_{A\sqcup B}\cdot x=0$. $\Box$\\ \\
Let $y$ be any of the generators $x,u_{A},v_{A},w_{A}$. As $\mbox{ad}(y)^{2}=0$ it follows that 
      $$(1+\mbox{ad}(y))^{2}=1+2\mbox{ad}(y)+\mbox{ad}(y)^{2}=1.$$
Thus $1+\mbox{ad}(y)$ is an involution in $\mbox{GL}(V^{*})$. \\ \\
{\bf Remark}. Notice that for any $A,B\subseteq {\mathbb N}$, the pairs $(\mbox{ad}(u_{A})$, $\mbox{ad}(u_{B}))$, $(\mbox{ad}(v_{A}),\mbox{ad}(v_{B}))$ and $(\mbox{ad}(w_{A}),\mbox{ad}(w_{B}))$ consist of elements that commute. Thus the subgroups 
        $${\mathcal U}=\langle 1+\mbox{ad}(u_{A}):\,A\subseteq {\mathbb N}\rangle,\,
				{\mathcal V}=\langle 1+\mbox{ad}(v_{A}):\,A\subseteq {\mathbb N}\rangle,\ 
				{\mathcal W}=\langle 1+\mbox{ad}(w_{A}):\,A\subseteq {\mathbb N}\rangle
				$$
are elementary abelian of countably infinite rank. We will be working with the group $G=\langle 1+\mbox{ad}(x), {\mathcal U},{\mathcal V},{\mathcal W}\rangle $. \\
\begin{lemm} The following commutator relations hold in $G$: \\
(a) $[1+\mbox{ad}(u_{A}),1+\mbox{ad}(v_{B})]=1+\mbox{ad}(u_{A\sqcup B})$. \\
(b) $[1+\mbox{ad}(v_{A}),1+\mbox{ad}(w_{B})]=1+\mbox{ad}(w_{A\sqcup B})$. \\
(c) $[1+\mbox{ad}(w_{A}),1+\mbox{ad}(u_{B})]=1+\mbox{ad}(v_{A\sqcup B})$. \\
(d) $[1+\mbox{ad}(u_{A}),1+\mbox{ad}(x)]=1.$ \\
(e) $[1+\mbox{ad}(v_{A}),1+\mbox{ad}(x)]=1.$ \\
(f) $[1+\mbox{ad}(w_{A}),1+\mbox{ad}(x)]=1+\mbox{ad}(u_{A})$. 
\end{lemm}
{\bf Proof}\  \ (a) We have
\begin{align*}
[1+\mbox{ad}(u_{A}),1+\mbox{ad}(v_{B})] & = (1+\mbox{ad}(u_{A}))\cdot (1+\mbox{ad}(v_{B}))\cdot 
                             (1+\mbox{ad}(u_{A}))\cdot (1+\mbox{ad}(v_{B})) \\
   & =  1+\mbox{ad}(u_{A})\mbox{ad}(v_{B})+\mbox{ad}(v_{B})
\mbox{ad}(u_{A}) \\
 & =  1+\mbox{ad}(u_{A}v_{B}) \\
& =  1+\mbox{ad}(u_{A\sqcup B}).
\end{align*}
(b) and (c) are proved similarly. For (f) we have
\begin{align*}
[1+\mbox{ad}(w_{A}),1+\mbox{ad}(x)] &=  (1+\mbox{ad}(w_{A}))\cdot (1+\mbox{ad}(x))\cdot 
                                        (1+\mbox{ad}(w_{A}))\cdot (1+\mbox{ad}(x)) \\
& =  1+\mbox{ad}(w_{A})\mbox{ad}(x)+\mbox{ad}(x)\mbox{ad}(w_{A})+\mbox{ad}(x)\mbox{ad}(w_{A})\mbox{ad}(x) \\
           & = 1+\mbox{ad}(w_{A}\cdot x) \\
           & = 1+\mbox{ad}(u_{A}).
\end{align*}
Here in the 2nd last equality, we have used Lemma 2.2. Parts (d) and (e) are proved similiarly. $\Box$ \\ \\
{\bf Remark}. Notice that as $V^{*}$ is locally nilpotent, it follows from Lemma 2.3 that $G$ is locally nilpotent. Next proposition clarifies further the 
structure of $G$. 
\begin{prop} We have $G=\langle 1+\mbox{ad}(x)\rangle {\mathcal U}{\mathcal V}{\mathcal W}$. Furthermore every
element $g\in G$ has a unique expression $g=(1+\mbox{ad}(x))^{\epsilon}rst$ with $\epsilon\in \{0,1\}$, $r\in {\mathcal U}$, $s\in {\mathcal V}$ and $t\in {\mathcal W}$.
\end{prop}
{\bf Proof}\ \ We first deal with the existence of such a decomposition. Suppose 
                  $$g=g_{0}(1+\mbox{ad}(x))g_{2}\cdots (1+\mbox{ad}(x))g_{n}$$
where $g_{0},\ldots ,g_{n}$ are products of elements of the form $1+\mbox{ad}(u_{A}),1+\mbox{ad}(v_{A})$
and $1+\mbox{ad}(w_{A})$. From Lemma 2.3 we know that $(1+\mbox{ad}(w_{A}))(1+\mbox{ad}(x))=
(1+\mbox{ad}(x))(1+\mbox{ad}(w_{A}))(1+\mbox{ad}(u_{A}))$ and $1+\mbox{ad}(x)$ commutes with all products of the form $1+\mbox{ad}(u_{A})$ and $1+\mbox{ad}(v_{A})$. We can thus collect the $(1+\mbox{ad}(x))$'s to the
left, starting with the leftmost occurrence. This may introduce more elements of the form $(1+\mbox{ad}(u_{A}))$ but no new $1+\mbox{ad}(x)$. We thus see that 
              $$g=(1+\mbox{ad}(x))^{n}g_{1}\cdots g_{m}$$
where $g_{i}$ is of the form $1+\mbox{ad}(u_{A})$, $1+\mbox{ad}(v_{A})$ or $1+\mbox{ad}(w_{A})$. This reduces our problem to the case when $g\in \langle {\mathcal U},{\mathcal V}, {\mathcal W}\rangle$. Suppose
       $$g=g_{0}(1+\mbox{ad}(u_{A_{1}}))g_{1}\cdots (1+\mbox{ad}(u_{A_{n}}))g_{n}$$ 
where $g_{0},\ldots, g_{n}$ are products of elements of the form $1+\mbox{ad}(v_{A})$ and 
$1+\mbox{ad}(w_{A})$. Suppose that the elements occurring in these products are $1+\mbox{ad}(v_{A_{n+1}}),
\ldots, 1+\mbox{ad}(v_{A_{n+l}}),1+\mbox{ad}(w_{A_{n+l+1}}),\ldots ,1+\mbox{ad}(w_{A_{m}})$. \\ \\
Using Lemma 2.3 we know that $(1+\mbox{ad}(v_{B}))(1+\mbox{ad}(u_{A}))=
(1+\mbox{ad}(u_{A}))(1+\mbox{ad}(v_{B}))(1+\mbox{ad}(u_{A\sqcup B}))$ and that $(1+\mbox{ad}(w_{B}))
(1+\mbox{ad}(u_{A}))=(1+\mbox{ad}(u_{A}))(1+\mbox{ad}(w_{B}))(1+\mbox{ad}(v_{A\sqcup B}))$. We can thus collect $1+\mbox{ad}(u_{A_{1}}),
\ldots ,1+\mbox{ad}(u_{A_{n}})$ to the left. In doing so we may introduce new terms of the form
$1+\mbox{ad}(u_{A})$, with $A$ of the form $A_{i_{1}}\sqcup \cdots \sqcup A_{i_{s}}$, and $s\geq 2$. This shows that 
      $$g=(1+\mbox{ad}(u_{A_{1}}))\cdots (1+\mbox{ad}(u_{A_{n}}))g_{1}\cdots g_{m}$$
where each $g_{j}$ is of the form $1+\mbox{ad}(v_{B}), 1+\mbox{ad}(w_{B})$ or $1+\mbox{ad}(u_{A})$, and  
$A$ is a modified union of at least $2$ sets from $\{A_{1},\ldots, A_{m}\}$. We can repeat this procedure, collecting all the new $(1+\mbox{ad}(u_{A}))$s. In doing so, we possibly introduce some new such elements but these will
then be with an $A$ that is a modified union of at least $3$ sets from $\{A_{1},\ldots ,A_{m}\}$. Continuing like this the procedure will
end after at most $m$ steps as every modified union of $m+1$ sets from $\{A_{1},\ldots ,A_{m}\}$ will be trivial. We have thus seen that $g=rh$ with $r\in {\mathcal U}$ and $h\in \langle {\mathcal V},{\mathcal W}\rangle$. We are now only left with the situation when $g\in \langle {\mathcal V},{\mathcal W}\rangle$. Suppose
       $$g=g_{0}(1+\mbox{ad}(v_{A_{1}}))g_{1}\cdots (1+\mbox{ad}(v_{A_{n}}))g_{n}$$
where $g_{0},g_{1},\ldots ,g_{n}$ are of the form $1+\mbox{ad}(w_{A})$. As $(1+\mbox{ad}(w_{B}))(1+\mbox{ad}(v_{A}))=(1+\mbox{ad}(v_{A}))(1+\mbox{ad}(w_{B}))(1+\mbox{ad}(w_{A\sqcup B}))$, we can now collect $1+\mbox{ad}(v_{A_{1}}), \ldots ,1+\mbox{ad}(v_{A_{n}})$ to the left and in doing so, all the new terms introduced will
be of the form $1+\mbox{ad}(w_{A})$. Thus $g=st$ with $s\in {\mathcal V}$ and $t\in {\mathcal W}$. This completes the existence part. We now want to show that such a decomposition is unique. Suppose 
\begin{eqnarray*}
    (1+\eta\, \mbox{ad}(x)) &  & (1+\tau\, \mbox{ad}(x)) \\
		  (1+\alpha_{1}\,\mbox{ad}(u_{A_{1}}))\cdots (1+\alpha_{r}\,\mbox{ad}(u_{A_{r}}))
		     &   &  (1+\beta_{1}\,\mbox{ad}(u_{A_{1}}))\cdots (1+\beta_{r}\,\mbox{ad}(u_{A_{r}}))\\
		(1+\gamma_{1}\,\mbox{ad}(v_{B_{1}}))\cdots (1+\gamma_{s}\,\mbox{ad}(v_{B_{s}})) & = &  (1+\delta_{1}\,\mbox{ad}(v_{B_{1}}))\cdots (1+\delta_{s}\,\mbox{ad}(v_{B_{s}})) \\
		(1+\epsilon_{1}\,\mbox{ad}(w_{C_{1}}))\cdots (1+\epsilon_{t}\,\mbox{ad}(w_{C_{t}})) &  & (1+\nu_{1}\,\mbox{ad}(w_{C_{1}}))\cdots (1+\nu_{t}\,\mbox{ad}(w_{C_{t}})).
\end{eqnarray*}
Applying both sides to $w_{\mathbb N}$ we get 
                $$w_{\mathbb N}+\eta\,u_{\mathbb N}=w_{\mathbb N}+\tau\,u_{\mathbb N}$$
from which we get $\eta=\tau$. Applying both sides to $x$ we get 
\begin{eqnarray*}
     x+\epsilon_{1}u_{C_{1}}+\cdots +\epsilon_{r}u_{C_{t}} & & x+\nu_{1}u_{C_{1}}+\cdots +\nu_{r}u_{C_{t}}\\
   \epsilon_{1}\epsilon_{2}v_{C_{1}\sqcup C_{2}}+\cdots +\epsilon_{t-1}\epsilon_{t}v_{C_{t-1}\sqcup C_{t}} & = & \nu_{1}\nu_{2}v_{C_{1}\sqcup C_{2}}+\cdots +\nu_{t-1}\nu_{t}v_{C_{t-1}\sqcup C_{t}} \\
 \epsilon_{1}\epsilon_{2}\epsilon_{3}w_{C_{1}\sqcup C_{2}\sqcup C_{3}}+\cdots  &  & 
	\nu_{1}\nu_{2}\nu_{3}w_{C_{1}\sqcup C_{2}\sqcup C_{3}}+ \cdots \\
	+
   \epsilon_{t-2}\epsilon_{t-1}\epsilon_{t}w_{C_{t-2}\sqcup C_{t-1}\sqcup C_{t}} & & +
   \nu_{t-2}\nu_{t-1}\nu_{t}w_{C_{t-2}\sqcup C_{t-1}\sqcup C_{t}} 
\end{eqnarray*}
from which we see that $\epsilon_{1}=\nu_{1},\ldots ,\epsilon_{t}=\nu_{t}$. Thus 
\begin{eqnarray*}
		  (1+\alpha_{1}\,\mbox{ad}(u_{A_{1}}))\cdots (1+\alpha_{r}\,\mbox{ad}(u_{A_{r}}))
		     &   &  (1+\beta_{1}\,\mbox{ad}(u_{A_{1}}))\cdots (1+\beta_{r}\,\mbox{ad}(u_{A_{r}}))\\
		(1+\gamma_{1}\,\mbox{ad}(v_{B_{1}})\cdots (1+\gamma_{s}\,\mbox{ad}(v_{B_{s}})) & = &  (1+\delta_{1}\,\mbox{ad}(v_{B_{1}})\cdots (1+\delta_{s}\,\mbox{ad}(v_{B_{s}})). 
\end{eqnarray*}
We can assume that $A_{j}\not\subseteq A_{i}$ and $B_{j}\not\subseteq B_{i}$ when $i<j$. Applying both sides
to $u_{{\mathbb N}\setminus B_{1}}$ gives 
      $$u_{{\mathbb N}\setminus B_{1}}+\gamma_{1}u_{\mathbb N}=u_{{\mathbb N}\setminus B_{1}}+\delta_{1}u_{\mathbb N}$$
from which we see that $\gamma_{1}=\delta_{1}$. Cancelling on both sides by $1+\gamma_{1}\,\mbox{ad}(v_{B_{1}})$ and then applying both sides to $u_{{\mathbb N}\setminus B_{2}}$ likewise gives $\gamma_{2}=\delta_{2}$. Continuing in this manner gives $\gamma_{1}=\delta_{1},\ldots ,\gamma_{s}=\delta_{s}$. We then have
$$(1+\alpha_{1}\,\mbox{ad}(u_{A_{1}}))\cdots (1+\alpha_{r}\,\mbox{ad}(u_{A_{r}}))=(1+\beta_{1}\,\mbox{ad}(u_{A_{1}}))\cdots (1+\beta_{r}\,\mbox{ad}(u_{A_{r}})).$$
A similar argument as before, applying both sides to $v_{{\mathbb N}\setminus A_{1}}, v_{{\mathbb B}\setminus A_{2}},\ldots $ gives likewise $\alpha_{1}=\beta_{1},\ldots ,\alpha_{r}=\beta{r}$. This finishes the
proof. $\Box$ \\ \\
We are now ready to prove the main result of this paper.
\begin{theo} The element $1+\mbox{ad}(x)$ is a left $3$-Engel element in $G$. However $\langle 1+\mbox{ad}(x)\rangle^{G}$ is not nilpotent. 
\end{theo}
{\bf Proof}\ \ Let $g=h(1+\mbox{ad}(w_{A_{1}}))\cdots (1+\mbox{ad}(w_{A_{n}}))$ be an arbitrary element 
in $G$ where $h\in \langle (1+\mbox{ad}(x))\rangle {\mathcal U}{\mathcal V}$. We want to show that 
$$[(1+\mbox{ad}(x))^{g},_{2}(1+\mbox{ad}(x))]=1.$$ Notice first that if $y\in V$ then 
\begin{eqnarray*}
     (1+\mbox{ad}(y))^{1+\footnotesize \mbox{ad}(w_{A})} & = & (1+\mbox{ad}(w_{A}))(1+\mbox{ad}(y))(1+\mbox{ad}(w_{A})) \\
                                         & = & 1+\mbox{ad}(y)+\mbox{ad}(yw_{A}).
\end{eqnarray*}
Notice that $(1+\mbox{ad}(x))^{g}=(1+\mbox{ad}(x))^{(1+\mbox{ad}(w_{A_{1}}))\cdots (1+\mbox{ad}(w_{A_{n}}))}$ and straighforward induction show that 
      $$(1+\mbox{ad}(x))^{g}=1+\mbox{ad}(y)$$
where 
 \begin{eqnarray*}
	y & = & x+\sum_{1\leq i\leq n}u_{A_{i}}+\sum_{1\leq i<j\leq n}v_{A_{i}\sqcup A_{j}}+
								\sum_{1\leq i<j<k\leq n}w_{A_{i}\sqcup A_{j}\sqcup A_{k}}.
\end{eqnarray*}
Since $\mbox{ad}(x)\mbox{ad}(y)\mbox{ad}(x)=0$, the commutator of $(1+\mbox{ad}(x))^{g}$ with $(1+\mbox{ad}(x))$ is
\begin{align*}
     (1+\mbox{ad}(y))(1+\mbox{ad}(x))(1+\mbox{ad}(y))(1+\mbox{ad}(x)) & = 
                1+\mbox{ad}(y)\mbox{ad}(x)+\mbox{ad}(x)\mbox{ad}(y)\\
              &   + \mbox{ad}(y)\mbox{ad}(x)\mbox{ad}(y).
\end{align*}
Then 
\begin{align*}
     [(1+\mbox{ad}(x))^{g},_{2}1+\mbox{ad}(x)] & = ((1+\mbox{ad}(y))(1+\mbox{ad}(x))^{4} \\ 
                          & = ((1+\mbox{ad}(y)\mbox{ad}(x)+\mbox{ad}(x)\mbox{ad}(y) \\ &+  \mbox{ad}(y)\mbox{ad}(x)\mbox{ad}(y))^{2}=1
\end{align*}
using again the fact that $\mbox{ad}(x)\mbox{ad}(y)\mbox{ad}(x)=0$. \\ \\
Then the normal closure of $1+\mbox{ad}(x)$ in $G$ is though not nilpotent as for $A_{i}=\{i\}$ we have 
  \begin{align*} [1+\mbox{ad}(w_{A_{1}}),1+\mbox{ad}(x),1+\mbox{ad}(w_{A_{2}}),1+\mbox{ad}(w_{A_{3}}),
& \ldots \\  \ldots, 1+\mbox{ad}(x),1+\mbox{ad}(w_{A_{2n}}),1+\mbox{ad}(w_{A_{2n+1}})]&=1+\mbox{ad}(w_{A}),
  \end{align*}
where $A=A_{1}\sqcup \ldots \sqcup A_{2n+1}=\{1,2,\ldots ,2n+1\}$. $\Box$
\mbox{}\\ \\
Our next aim is to show however that if we take any $r$ conjugates $(1+\mbox{ad}(x))^{g_{1}},\ldots,(1+\mbox{ad}(x))^{g_{r}}$ of $(1+\mbox{ad}(x))$ in $G$, they generate a nilpotent subgroup of $r$-bounded class that grows linearly with $r$. \\ \\
We first work in a more general setting. For each $e\in E$ and $\emptyset \not =A\subseteq {\mathbb N}$, let $e(A)\in \mbox{End}(V^{*})$ where
                         $$u_{B}e(A)=(ue)_{B\sqcup A}.$$
Then let $E^{*}=\langle \mbox{ad}(x), e(A):\,e\in E\mbox{ and }\emptyset \not =A\subseteq {\mathbb N}\rangle$. As $V^{*}$ is locally nilpotent, one sees readily that
the elements of $E^{*}$ are nilpotent and thus $1+E^{*}$ is a subgroup of $\mbox{End}(V^{*})$. We are going to see that $1+E^{*}$ is of finite 
exponent. \\ 
{\bf Remark} Notice that $\mbox{ad}(u_{\mathbb N}) = \mbox{ad}(v_{\mathbb N}) = 0.$
\begin{lemm} The elements $\mbox{ad}(w_{\mathbb N})$ and $\{e_{i}(A):\,1\leq i\leq 12,\,\emptyset \not =A\subset {\mathbb N}\} \cup \{e_1({\mathbb N}), e_2({\mathbb N}), e_3({\mathbb N})\}$ form a basis for $E^{*}$.
\end{lemm}
{\bf Proof}\ \ One sees that these elements span $E^{*}$ as a vector space in a similar way as in the proof of Lemma 2.1. We thus skip over the details
and only show that these elements are linearly independent. Suppose 
                  $$\epsilon\mbox{ad}(x) +\sum_{i=1}^{12}\sum_{A}\epsilon^{i}_{A}e_{i}(A)=0,$$
where only finitely many of the coefficients $\epsilon, \epsilon^{i}_{A}$ are non-zero. Denote the left hand
side by $T$. Then 
    $$0=xT=\sum_{A}\epsilon^{1}_{A}u_{A}+\sum_{A}\epsilon^{2}_{A}v_{A}+\sum_{A}\epsilon^{3}_{A}w_{A}$$
implying that $\epsilon^{i}_{A}=0$ for all $A$ and $i=1,2,3$. Then 
    $$0=u_{{\mathbb N}\setminus A}T=\sum_{B\subseteq A}\epsilon^{4}_{B}u_{({\mathbb N}\setminus A)\sqcup B}+
		\sum_{B\subseteq A}\epsilon^{5}_{B}v_{({\mathbb N}\setminus A)\sqcup B}+
		\sum_{B\subseteq A}\epsilon^{6}_{B}w_{({\mathbb N}\setminus A)\sqcup B}.$$
In particular $\epsilon^{i}_{A}=0$ for all $A \neq \mathbb{N}$ and $i=4,5,6$. We continue in a similar way. Next
  $$0=v_{{\mathbb N}\setminus A}T=\sum_{B\subseteq A}\epsilon^{7}_{B}u_{({\mathbb N}\setminus A)\sqcup B}+
		\sum_{B\subseteq A}\epsilon^{8}_{B}v_{({\mathbb N}\setminus A)\sqcup B}+
		\sum_{B\subseteq A}\epsilon^{9}_{B}w_{({\mathbb N}\setminus A)\sqcup B}$$ 
that shows that $\epsilon^{i}_{A}=0$ for all $A$ and $i=7,8,9$. Finally 
                $$0=w_{\mathbb N}T=\epsilon u_{\mathbb N}$$
giving $\epsilon=0$ and 
$$0=w_{{\mathbb N}\setminus A}T=\sum_{B\subseteq A}\epsilon^{10}_{B}u_{({\mathbb N}\setminus A)\sqcup B}+
		\sum_{B\subseteq A}\epsilon^{11}_{B}v_{({\mathbb N}\setminus A)\sqcup B}+
		\sum_{B\subseteq A}\epsilon^{12}_{B}w_{({\mathbb N}\setminus A)\sqcup B}$$ 
and $\epsilon^{i}_{A}=0$ for all $A$and $i=10,11,12$. This finishes the proof. $\Box$
\begin{coro} We have $(1+E^{*})^{32}=1$.
\end{coro}
{\bf Proof}\ \ Let $\bar{E}$ be the subalgebra of $E^{*}$ generated by all $e_{i}(A)$ where $1\leq i\leq 12$
and $\emptyset \not =A\subseteq {\mathbb N}$. Let
$f=\mbox{ad}(x)+e\in E$ where $e\in \bar{E}$. Then $f^{2}=\mbox{ad}(x)^{2}+e^{2}+(e\mbox{ad}(x)+\mbox{ad}(x)e)=
e^{2}+(e\mbox{ad}(x)-\mbox{ad}(x)e)$. It is straightforward to see that $e\mbox{ad}(x)-\mbox{ad}(x)e\in \bar{E}$ ($\bar{E}$ is an ideal in the Lie algebra $E^{*}$) and thus $f^{2}\in \bar{E}$. 
It thus suffices to show that $\bar{E}^{16}=0$, as then it will follow that $(E^{*})^{32}=0$ and therefore $(1+e)^{32}=1+e^{32}=1$ for
all $e\in E^{*}$. \\ \\
Let $e=y_{1}+\cdots +y_{m}$ be any element in $\bar{E}$ where $y_{1},\ldots, y_{m}$ belong to the basis $\{e_{i}(A):\,1\leq i\leq 12,\,\emptyset \not =A\subset {\mathbb N}\} \cup \{e_1({\mathbb N}), e_2({\mathbb N}), e_3({\mathbb N})\}$ for $\bar{E}$ given in Lemma 2.6. As any product with a repeated term is $0$ we see that $e^{16}$ is a sum of terms of the form
\begin{equation}
            \sum_{\sigma\in S_{16}}f_{\sigma(1)}\cdots f_{\sigma(16)}
\end{equation}
with $f_{1},\ldots ,f_{16}\in \{y_{1},\ldots ,y_{m}\}$. As $16>12$ some two of $f_{1},\ldots, f_{16}$ must
be of the same type.  Without loss of generality we can assume that these are $f_{15}=e_{i}(A)$ and
$f_{16}=e_{i}(B)$. Notice that the sum (1) splits naturally into $16!/2$ sums of pairs
 \begin{eqnarray*}
\sum_{\sigma\in S_{16}}f_{\sigma(1)}\cdots f_{\sigma(16)} & = & \sum_{\sigma\in S_{14}}(e_{i}(A)e_{i}(B)+e_{i}(B)e_{i}(A))f_{\sigma(1)}\cdots f_{\sigma(14)}+\cdots \\ 
   & & + \sum_{\sigma\in S_{14}}f_{\sigma(1)}\cdots f_{\sigma(14)}(e_{i}(A)e_{i}(B)+e_{i}(B)e_{i}(A)),
\end{eqnarray*}
one for each of the ${16\choose 2}$ positions of the pair $(e_{i}(A),e_{i}(B))$ within the product. But for each such choice of positions the two elements in the pair have the same value and as the characteristic is $2$, the sum of each pair is $0$. Thus the sum in (1) is zero and we have shown that $\bar{E}^{16}=0$. $\Box$
\begin{prop} Any $r$-generator subgroup of $1+E^{*}$ is nilpotent of $r$-bounded class. 
\end{prop}
{\bf Proof} From Corollary 2.7 we know that $1+E^{*}$ is of bounded exponent. The result thus follows from Zel'manov's solution to the Restricted Burnside Problem. $\Box$ \\ \\ 
Despite the fact that the normal closure of $1+\mbox{ad}(x)$ in $G$ is not nilpotent, it turns out that the nilpotency class of the subgroup
generated by any $k$ conjugates grows linearly with respect to $k$. In order to see this we first introduce some more notation. Let $A_{1}, A_{2},
\ldots ,A_{k}$ be any $k$ subsets of ${\mathbb N}$. For each $k$-tuple $(i1,i2,\ldots ,ir)$ of non-negative integers and each $e\in E$ we let 
                 $$e^{(i1,\ldots ,ir)}=\sum_{\tiny \begin{array}{l}B_{1}\subseteq A_{1}\\
								                               |B|=i1																							\end{array}}
																							\ldots 
																							\sum_{\tiny \begin{array}{l}B_{k}\subseteq A_{k}\\
								                               |B|=ik																					\end{array}}e(B_{1}\sqcup B_{2}\sqcup \cdots \sqcup B_{k}).$$
Notice that        
         $$e^{(i1,\ldots ,ir)}f^{(j1,\ldots ,jr)}={\tiny i1+j1\choose i1}\cdots {\tiny ir+jr\choose ir}\normalsize (ef)^{(i1+j1,\ldots ,ir+jr)}. $$
Now notice that ${3+i\choose 3}$ is even for $i=1,2,3$ and the same is true for ${2+2\choose 2}$ and ${1+1\choose 1}$. However ${2+1\choose 2}$ is odd. From this it follows that 
                $$Q=\langle \mbox{ad}(x),e^{(i1,\ldots ,ir)}:\,e\in E,\, 0\leq i1,\ldots ,ir\leq 3,\ i1+\cdots+ir\geq 1\rangle$$
is a subalgebra of $E^{*}$. A nonzero product in $Q$ can have at most $2r$ elements of the form $e^{(i1,\ldots ,ir)}$ and as $\mbox{ad}(x)^{2}=0$ we could then have at most $1+2r$ occurrences of $\mbox{ad}(x)$ in a non-zero product. Thus $Q^{4r+2}=0$. \\ \\
Now take some $r$ conjugates of $(1+\mbox{ad}(x))$ in $G$. Recall that each conjugate is of the from 
$(1+\mbox{ad}(x))^{(1+\footnotesize{\mbox{ad}(w_{C_{1}}})\cdots (1+\footnotesize{\mbox{ad}(w_{C_{j}}}))}$. For ease of notation we will assume that each $C_{k}$ is a singleton set. The following argument also works for the more general case. Let 
     $$A_{1}=\{1,\ldots ,k_{1}\},\ A_{2}=\{k_{1}+1,\ldots ,k_{2}\},\ldots , A_{r}=\{k_{r-1}+1,\ldots,k_{r}\}.$$																
Then we have seen (see the proof of Theorem 2.5) that 
\begin{eqnarray*}
        (1+\mbox{ad}(x))^{(1+\mbox{ad}(w_{1}))\cdots (1+\mbox{ad}(w_{k_{1}}))} & = &
				   1+\mbox{ad}(x)+e_{7}^{(1,0,\ldots ,0)}+e_{4}^{(2,0,\ldots, 0)} + e_{1}^{(3,0,\ldots ,0)}\\
                                     & \vdots & \\
        (1+\mbox{ad}(x))^{(1+\mbox{ad}(w_{k_{r-1}+1}))\cdots (1+\mbox{ad}(w_{k_{r}}))}																		& = & 1+\mbox{ad}(x)+e_{7}^{(0,\ldots ,0,1)}+
				e_{4}^{(0,\ldots ,0,2)}+e_{1}^{(0,\ldots ,0,3)}.
\end{eqnarray*}
In other words the $r$ conjugates are all in $1+Q$. Hence if $H$ is the subgroup of $\mbox{GL}(V^{*})$ generated by the $r$ conjugates then 
                               $$\gamma_{4r+2}(H)\leq \gamma_{4r+2}(1+Q)\leq 1+Q^{4r+2}=1.$$
We have thus proved the following.
\begin{prop}Let $(1+\mbox{ad}(x))^{g_{1}},\ldots ,(1+\mbox{ad}(x))^{g_{r}}$ be any $r$ conjugates of $1+\mbox{ad}(x)$ in $G$. Then the group generated by these conjugates is nilpotent of class at most $4r+2$. 
\end{prop}

\subsubsection*{Acknowledgments}
The first author is partially supported by the ``National Group for Algebraic and Geometric Structures, and their Applications'' (GNSAGA \-- INdAM). Also, she would like to thank the Department of Mathematics at the University of Bath for its excellent hospitality while this research was conducted.

The second and the third author acknowledge the EPSRC (grant number 1652316) for their support.

\end{document}